\documentclass[graybox]{svmult}
\usepackage[OT2,T1]{fontenc}
\usepackage[utf8]{inputenc}

\usepackage{mathptmx}       
\usepackage{helvet}         
\usepackage{courier}        
\usepackage{type1cm}        
\usepackage{makeidx}         
\usepackage{graphicx}        
\graphicspath{{Figures/}}
\usepackage{multicol}        
\usepackage[bottom]{footmisc}

\makeindex             

\usepackage{filecontents}
\usepackage[colorlinks=true]{hyperref}
  \hypersetup{
    pdftitle=On Multilateral Hierarchical Dynamic Decisions, 
    pdfauthor=K. Szajowski, 
    colorlinks,
    urlcolor=blue,
    filecolor=magenta,
    citecolor=green, 
    linkbordercolor={1 1 1}, 
    citebordercolor={1 1 1},  
    urlbordercolor={ 1 1 1}  
  }
\usepackage{amsmath}
\usepackage{amssymb}

\RequirePackage[numbers]{natbib}

\newcommand*{\doi}[1]{\href{http://dx.doi.org/#1}{doi: #1}}
\newcommand*{\MR}[1]{\href{http://www.ams.org/mathscinet-getitem?mr=#1&return=pdf}{MR #1}}
\newcommand*{\ZBL}[1]{\href{http://www.zentralblatt-math.org/zmath/en/advanced/?q=an:#1&format=complete}{Zbl #1}}

\def\one{{\mathbb I}}
\def\bbN{{\mathbb N}}

\def\bE{\mathbf{E}}
\def\bP{\mathbf{P}}

\def\cA{{\mathcal A}}
\def\cB{{\mathcal B}}
\def\cF{{\mathcal F}}
\def\cG{{\mathcal G}}
\def\cH{{\mathcal H}}
\def\cP{{\mathcal P}}

\def\cS{{\mathcal S}}
\def\cT{{\mathcal T}}

\def\fM{{\mathfrak M}}

\def\tX{\tilde X}

\newcommand\myatop[2]{\genfrac{}{}{0pt}{}{#1}{#2}}
\begin{document}

\title*{On Multilateral Hierarchical Dynamic Decisions}
\titlerunning{Multilateral Decision models} 
\author{Krzysztof Szajowski}
\authorrunning{K. Szajowski} 
\institute{Krzysztof Szajowski \at Wrocław University of Science and Technology, Wybrzeże Wyspiańskiego 27, 50-370 Wrocław, Poland \email{Krzysztof.Szajowski@pwr.edu.pl}
}
%
%
\maketitle
\abstract*{Many decision problems in economics, information technology and industry can be transformed to an optimal stopping of adapted random vectors with some utility function over the set of Markov times with respect to filtration build by the decision maker's knowledge. The optimal stopping problem formulation is to find a stopping time which maximizes the expected value of the accepted (stopped) random vector's utility.\\
There are natural extensions of optimal stopping problem to stopping games-the problem of stopping random vectors by two or more decision makers. Various approaches dependent on the information scheme and the aims of the agents in a considered model. This report unifies a group of non-cooperative stopping game models with forced cooperation by the role of the agents, their aims and aspirations (v. Assaf \& Samuel-Cahn(1998), Szajowski \& Yasuda(1995)) or extensions of the strategy sets (v. Ramsey \& Szajowski(2008)).}

\abstract{Many decision problems in economics, information technology and industry can be transformed to an optimal stopping of adapted random vectors with some utility function over the set of Markov times with respect to filtration build by the decision maker's knowledge. The optimal stopping problem formulation is to find a stopping time which maximizes the expected value of the accepted (stopped) random vector's utility.\newline\indent
There are natural extensions of optimal stopping problem to stopping games-the problem of stopping random vectors by two or more decision makers. Various approaches dependent on the information scheme and the aims of the agents in a considered model. This report unifies a group of non-cooperative stopping game models with forced cooperation by the role of the agents, their aims and aspirations (v. \citeauthor{asssam98:coop}(\citeyear{asssam98:coop}), \citeauthor{szayas95:voting}(\citeyear{szayas95:voting})) or extensions of the strategy sets (v. \citeauthor{ramsza08:correlated}(\citeyear{ramsza08:correlated})).
\vskip3pt \noindent
}

\section{\label{sec:1}Introduction}
The subject of the analysis is the problem of making collective decisions by the team of agents in which the position (significance) of the members is not equal. An object that is subject to management generates a signal that changes over time. Agents deal with capturing signals. Everyone can capture and save one of them, and its value is relative, determined by the function that takes into account the results of all decisions. Both the ability to observe signals and their capture determines the rank of agents who compete in this process. It is also possible that unequalizes of the decision makers is a consequence of social agreement or policy (v. \citeauthor{FenXia2000:Revenue}(\citeyear{FenXia2000:Revenue})).

Earlier studies  of such issues (cf. \citeauthor{Fer1992:2person}(\citeyear{Fer1992:2person}), \citeauthor{Sza1993:2S2DM}(\citeyear{Sza1993:2S2DM}), \citeauthor{RamCie2009:Cooperative}(\citeyear{RamCie2009:Cooperative}), \citeauthor{DorManPon2009:OS}(\citeyear{DorManPon2009:OS}), \citeauthor{KraFer2013:Elfving}(\citeyear{KraFer2013:Elfving}), \citeauthor{Fer2016:SGSakaguchi}(\citeyear{Fer2016:SGSakaguchi})) showed their complexity, and detailed models of the analyzed cases a way to overcome difficulties in modeling and setting goals with the help of created models. The basic difficulty, except for cases when the decision is made by one agent, consists in determining the goals of the team, which can not always be determined so that the task can be reduced to the optimization of the objective function as the result of scalarisation. Most often, individual agents are to achieve an individual goal, but without the destabilization of the team. When modeling such a case, one should remember about establishing the rational goal of the agents in connection with the existence of the team (v. \citeauthor{DieVen2008:Aspiration}(\citeyear{DieVen2008:Aspiration})). In the considerations of this study, we use methods of game theory with a finite number of players. However, the classic model of the antagonistic game is not the best example of progress. The team has interactions of agents resulting even from the hierarchy of access to information and the order in which decisions are made. The proposed overcoming of this difficulty consists in the appropriate construction of strategy sets and the payment function of players so that, taking into account the interactions, construct a multi-player game in which players have sets of acceptable strategies chosen regardless of the decisions of other players. Due to the sequential nature of the decision-making process, this player's decision-making independence is at the time of making it, but it is conditioned by the team's existing decision-making process.

Due to the fact that the goal of each agent, aspiration assessment by defining a withdrawal function, is to accept the most important signal from its point of view, the result of modeling is the task of repeatedly stopping the sequence of random vectors. In fact rating aspirations by defining the functions of payment is the one of the preliminary work on the mathematical modeling of management problem.  Taking this into account, it should be mentioned here that this task was first put forward by ~\citeauthor{hagg1967:2stop}(\citeyear{hagg1967:2stop}), although \citeauthor{dyn69:game}'s(\citeyear{dyn69:game}) considerations can also be included in this category. Despite the undoubtedly interesting implications of such a model in applications, the subject has not been explored too much in its most general formulation, at least it has not been referred to. We will try to point out considerations that support such implicit modeling\footnote{The stopping games as the sepcial case of the stochastic game has been presented by \citeauthor{JasNow2016:Non-zero}(\citeyear{JasNow2016:Non-zero})}. 

In the game models applied to business decisions there are important models formulated and investigated by economist \citeauthor{Sta1934:Marktform}(\citeyear{Sta1934:Marktform})\footnote{This is his habilitation (see also the dissertation \cite{Sta1932:Kostentheorie}),  translated recently to English and published by Springer \citeauthor{Sta2011:enMarket}(\citeyear{Sta2011:enMarket}) (v. \cite{Etro2013:ReviewStackelberg} for the review of the edition).}. Formulation of the game related to the secretary problem by Fushimi~\cite{fus1981:competitive} with restricted set of strategies, namely threshold stopping times, opened research on the stopping games with leader by Szajowski~(see papers~\citeauthor{RavSza1992:Dynkin}(\citeyear{RavSza1992:Dynkin}), ~\citeauthor{Sza1992:Priority}(\citeyear{Sza1992:Priority})). Similar games are subject of the research by \citeauthor{EnnFer1987:Priorities}(\citeyear{EnnFer1987:Priorities}), \citeauthor{RadSza1988:sequential}(\citeyear{RadSza1988:sequential}). The extension of idea of Stackelberg was assumption that the lider is not fixed but the priority to the player is assigned randomly. Such version of the stopping game is investigated in \citeauthor{Sza1993:ZOR}(\citeyear{Sza1993:ZOR}) \citeauthor{Sza1993:2S2DM}(\citeyear{Sza1993:2S2DM}) and \citeauthor{Sza1995:SIAM}(\citeyear{Sza1995:SIAM}). The Nash equilibria are obtained in the set of randomized strategies (cf. \citeauthor{neuramsza02:Randomized}(\citeyear{neuramsza02:Randomized}), \citeauthor{neuporsza94:note}(\citeyear{neuporsza94:note})). 

Two or multi-person process stoppages, originally formulated by \citeauthor{dyn69:game}(\citeyear{dyn69:game})\footnote{See also models created by \citeauthor{McK1965:Free}(\citeyear{McK1965:Free}), \citeauthor{Kif1969:Games}(\citeyear{Kif1969:Games}).}, met with more interest and research on multi-player games with stopping moments as players' strategies are quite well described in the literature. Both for random sequences and for certain classes of processes with continuous time. We will use this achievement in our deliberations. 

In the following chapters, \ref{sec:2}-\ref{sec:4}, we will discuss hierarchical diagrams in multi-person decision problems and their reduction to an antagonistic game. We will use the lattice properties of the stopping moments and we will obtain an equilibrium point in the problems under consideration based on the fixed point theorem for the game on the complete grating.\\[-4ex]

\section{\label{sec:2}Decision makers' hierarchy in multi-choice problem} 
Let us consider $N$ agents multiple-choice decision model on observation of stochastic sequence. The decision makers (DMs) are trying to choose the most profitable state based on sequential observation. In the case when more than one player would like to accept the state there are priority system which choose the beneficiary and the other players have right to  observe further states of the process trying to get their winning state. 

The aims of the agents are defined by the pay-offs function of the them. The rationality is subject of arbitrary decision when the mathematical model is formulated and should emphasize the requirement of the agents. One of the popular way is transformation of such multilateral problem to a non-zero-sum game. When there are two DMs it could be also zero-sum stopping game.\\[-4ex]

\subsection{Zero-sum Dynkin's Game}
The originally~\citeauthor{dyn69:game}(\citeyear{dyn69:game}) has formulated the following optimization problem. Two players observe a realization of two real-valued processes $(X_n)$ and $(R_n)$. Player~1 can stop whenever $X_n\geq 0$, and player~2 can stop whenever $X_n < 0$. At the first stage $\tau$ in which one of the players stops, player~2 pays player~1 the amount $R_\tau$ and the process terminates. If no player ever stops, player~2 does not pay anything. 

A \emph{strategy} of player~1 is a stopping time $\tau$ that satisfies $\{\tau = n\} \subset \{X_n\geq 0\}$ for every $n\geq 0$. A strategy $\sigma$ of player~2 is defined analogously. The termination stage is simply $\nu= \min\{\tau,\sigma\}$. For a given pair $(\tau,\sigma)$ of strategies, denote by $K(\tau,\sigma) = \bE\one_{\{\nu<\infty\}}R_\nu$ the expected payoff to player~1.

\citeauthor{dyn69:game}(\citeyear{dyn69:game}) proved that if $\sup_{n\geq 0} |R_n|\in L_1$ then this problem has a value $v$ i.e. 
\begin{equation*}
v = \sup_\tau\inf_\sigma K(\tau,\sigma)=\inf_\sigma\sup_\tau K(\tau,\sigma)
\end{equation*}
\vspace{-4ex}
\subsection{Non-zero sum stopping game}
Basement process under which the game is formulated can be defined as follows. Let $(X_n,\cF_n,\bP_x)_{n=0}^T$ be a homogeneous Markov process defined on a probability space $(\Omega,\cF,\bP)$ with state space $(\E,\cB)$.  At each moment $n =1,2,...,T$, $T\in \tilde{\bbN}=\bbN\cup\{\infty\}$, the decision makers (henceforth called players) are able to observe the consecutive states of Markov process sequentially.  There are $N$ players. Each player has his own utility function $g_i:\E^N\rightarrow \Re$, $i=1,2,\ldots,N$, dependent on his own and others choices of state the Markov process. At moment $n$ each decides separately whether to accept or reject the realization $x_n$ of $X_n$.   We assume the functions $g_i$ are measurable and bounded.

\begin{itemize}
\item Let $\cT^i$ be the set of pure strategies for $i$th player, the stopping times with respect to the filtration $(\cF_n^i)_{n=1}^T$, $i=1,2,\ldots,N$. Each player has his own sequence of $\sigma$-fields $(\cF_n^i)_{n=1}^T$ (the  available information).
\item The randomize extension of $\cT_i$ can be constructed as follows (see \citeauthor{yas85:randomize}(\citeyear{yas85:randomize}), \citeauthor{shmsol04:nonzero}(\citeyear{shmsol04:nonzero})). 
Let $(A_n^i)_{n=1}^T$, $i=1,2,\ldots,N$, be i.i.d.r.v. from the uniform distribution on $[0,1]$ and independent of the Markov process $(X_n,\cF_n,\bP_x)_{n=0}^T$. Let $\cH_n^i$ be the $\sigma$-field generated by $\cF_n^i$ and  $\{(A_s^i)_{s=1}^n \}$. A randomized Markov time $\tau(p^i)$ for strategy $p^i = (p_n^i)\in \cP^{T,i}\in\fM^T_i$ of the $i$th player is $\tau(p^i) = \inf\{T\geq n\geq 1: A_n^i\leq p_n^i\}$. 
\end{itemize}
\vspace{-4ex}

Clearly, if each $p_n^i$ is either zero or one, then the strategy is pure and $\tau(p^i)$ is in fact an $\{\cF_n^i\}$- Markov time. In particular an $\{\cF_n^i\}$- Markov time $\tau_i$ corresponds to the strategy $p^i = (p_n^i)$ with $p_n^i = \one_{\{\tau_i = n\}}$, where $\one_A$ is the indicator function for the set $A$. 

Two concepts are take into account in this investigation. It can be compared with real investments and investment on the financial market. In real investment the choice of state is not reversible and sharable. In the financial market the choice of state by many players can be split of profit to all of them according some rules. Here, it is separately considered models of payoffs definition.

The payoff functions should be adequate to the information which players have and their decision. The player who do not use his information should be penalize.  
\begin{itemize}
\item Let the players choose the strategies $\tau_i\in\cT^i$, $i=1,2,\ldots,N$. The payoff of the $i$th player is $G_i(\tau_1,\tau_2,\ldots,\tau_N)=g_i(X_{\tau_1},X_{\tau_2},\ldots,X_{\tau_N})$.
\item If the $i$th player control the $i$th component of the process, than the function $G_i(i_1,i_2,\ldots,i_n)=h_i(X_{i_1},X_{i_2},\ldots,X_{i_n})$ forms the random field. Such structure of payoffs has been considered by \citeauthor{mam87:monotone}(\citeyear{mam87:monotone}). Under additional assumptions concerning monotonicity of incremental benefits of players Mamer has proved the existence of Nash equilibrium for two player non-zero sum game. 
\item Let $\cG_n=\sigma(\cF_n^1\cup\cF_n^2\cup\ldots\cup\cF_n^N)$ and $\cT$ be the set of stopping with respect to $(\cG_n)_{n=1}^T$. For a given choice of strategies by players the effective stopping time $\nu=\psi(\tau_1,\tau_2,\ldots,\tau_N)$ and $G_i(\tau_1,\tau_2,\ldots,\tau_N)=g_i(X_\nu)$. In some models the process $X_n$ can be multidimensional and the payoff of $i$th player is the $i$th component of the vector $X_n$. 
\end{itemize}

\begin{definition}[Nash equilibrium]
The strategies $\tau_1^\star,\tau_2^\star,\ldots,\tau_N^\star$ are equilibrium in stopping game if for every player $i$ 
\begin{equation}
\bE_x G_i(\tau_1^\star,\tau_2^\star,\ldots,\tau_N^\star)\geq \bE_x G_i(\tau_1^\star,\tau_2^\star,\ldots,\tau_i,\ldots,\tau_N^\star).
\end{equation}
\end{definition}
\vspace{-4ex}
     
\subsection{Rights assignment models}
However, there are different systems of rights to collect information about underlined process and priority in acceptance the states of the process. The various structures of decision process can have influence the knowledge of the players about the process which determine the pay-offs of the players. It is assumed that the priority decide about the investigation of the process and decision of the state acceptance.
The details of the model, which should be precise are listed here.
\begin{enumerate}
\item The priority of the players can be defined before the game (in deterministic or random way) or it is dynamically managed in the play.
\item The priority of the players is decided after the collection of knowledge about the item by all players.
\item The random assignment of the rights can run before observation of each item and the accepted observation is not known to players with lowest priority. It makes that after the first acceptance some players are better informed than the others.  
\begin{enumerate}
\item The information about accepted state is known to all players.
\item The information is hidden to the players who do not accepted the item. 
\end{enumerate}
\end{enumerate}
The topics which are analyzed could be pointed out as follows:
 \begin{enumerate}
  \item Dynkin's game; 
  \item The fix and dynamic priority of the players: deterministic and random;
  \end{enumerate}  
\vspace{-4ex}  \nopagebreak
 
\section{\label{sec:3}The fix and dynamic priority of the players.}
\subsection{Deterministic priority}\subsubsection{Static.} Among various methods of privilages for the players one of the simples is permutation of players' indices (rang). Let us propose a model of assignments the priority (rang) to the players as follows. 
In non-zero two person Dynkin's game the role of an \emph{arbiter} was given to the random process ${X_n}$. The simplest model can assume that the players are ordered before the play to avoid the conflict in assignment of presented sequentially states. At each moment the successive state of the process is presented to the players, they decide to stop and accept the state or continue observation. The state is given to the players with highest rang (we adopt here the convention that the player with rang $1$ has the highest priority). In this case each stopping decision reduce the number of players in a game. It leads to recursive algorithm of construction the game value and in a consequence to determining the equilibrium (see \citeauthor{nowsza94:stochastic}(\citeyear{nowsza94:stochastic}), \citeauthor{sak95:review}(\citeyear{sak95:review}) for review of such models investigation).   

The players decision and their priorities define an effective stopping time for player $i$ in the following way.
\begin{itemize}
\item Let $P=\{1,2,\ldots,N\}$ be the set of players and $\pi$ a permutation of $P$. It determines the priority $\pi(i)$ of player $i$.  
\end{itemize}
The considered model can be extended to fix deterministic priority. The effective stopping time for player $i$ in this case one can get as follows.
\begin{itemize}
\item Let $(p_n^i)_{n=1}^T$ be the pure stopping strategy. If it is randomized stopping time we can find pure stopping time with respect to an extended filtration. The effective stopping strategy of the player $i$ is following:
\begin{equation}
\tau_i(\vec(p))=\inf\{k\geq 1:p_k^i\prod_{j=1}^N(1-p_k^j)\one_{\{j:\pi(j)<\pi(i)\}}=1\},
\end{equation}
where $\vec{p}=(p^1,p^2,\ldots,p^N)$ and each $p^i=(p^i_n)_{n=1}^T$ is adapted to the filtration $(\cF_n^i)_{n=1}^T$. The effective stopping time of the player $i$ is the stopping time with respect to the filtration $\tilde{\cF_n^i}=\sigma\{\cF_n^i,\{(p_k^j)_{k=1,\{j: \pi(j)<\pi(i)\}}^n\}\}$. 
\item The above construction of effective stopping time assures that each player will stop at different moment. It translates the problem of fixed priority optimization problem to the ordinary stopping game with payoffs $G_i(\tau_1,\tau_2,\ldots,\tau_N)=g_i(X_{\tau_1},X_{\tau_2},\ldots,X_{\tau_N})$. 
\end{itemize}
\subsubsection{Dynamic}
In this case the effective stopping time for player $i$ is obtained from parameters of the model similarly.
\begin{itemize}
\item Let $(p_n^i)_{n=1}^T$ be the pure stopping strategy. If it is randomized stopping time we can find pure stopping time with respect to an extended filtration. The effective stopping strategy of the player $i$ is following:
\begin{equation}
\tau_i(\vec{p})=\inf\{k\geq 1:p_k^i\prod_{j=1}^N(1-p_k^j)\one_{\{j:\pi_k(j)<\pi_k(i)\}}=1\},
\end{equation}
where $\vec{p}=(p^1,p^2,\ldots,p^N)$ and each $p^i=(p^i_n)_{n=1}^T$ is adapted to the filtration $(\cF_n^i)_{n=1}^T$. The effective stopping time of the player $i$ is the stopping time with respect to the filtration $\tilde{\cF_n^i}=\sigma\{\cF_n^i,\{(p_k^j)_{k=1,\{j: \pi_k(j)<\pi_k(i)\}}^n\}\}$. 
\item The above construction of effective stopping time assures that each player will stop at different moment. It translates the problem of fixed priority optimization problem to the ordinary stopping game with payoffs $G_i(\tau_1,\tau_2,\ldots,\tau_N)=g_i(X_{\tau_1},X_{\tau_2},\ldots,X_{\tau_N})$. 
\end{itemize}
\subsection{The random priority of the players}\subsubsection{Static(fixed) and dynamic}
The random permutation of the players' can be model of the random fix priority when before the play the assignment of priority is based on the random permutation. The fixed permutation is valid for one turn of the game. The effective stopping time for player $i$ has the following construction in this case.
\begin{itemize}
\item It is still fixed permutation of the player but its choice is random. The drawing of the permutation $\Pi$ is done onces for each play. Let $(p_n^i)_{n=1}^T$ be the pure stopping strategy. If it is randomized stopping time we can find pure stopping time with respect to an extended filtration. The effective stopping strategy of the player $i$ is following:
\begin{equation}
\tau_i(\vec{p})=\inf\{k\geq 1:p_k^i\prod_{j=1}^N(1-p_k^j)\one_{\{j:\Pi(j)<\Pi(i)\}}=1\},
\end{equation}
with rest of denotations the same as in the previous section, i.e. where $\vec{p}=(p^1,p^2,\ldots,p^N)$ and each $p^i=(p^i_n)_{n=1}^T$ is adapted to the filtration $(\cF_n^i)_{n=1}^T$. The effective stopping time of the player $i$ is the stopping time with respect to the filtration $\tilde{\cF_n^i}=\sigma\{\cF_n^i,\Pi,\{(p_k^j)_{k=1,\{j:\Pi(j)<\Pi(i)\}}^n\}\}$. 
\item The above construction of effective stopping time assures that each player will stop at different moment. It translates the problem of fixed priority optimization problem to the ordinary stopping game with payoffs $G_i(\tau_1,\tau_2,\ldots,\tau_N)=g_i(X_{\tau_1},X_{\tau_2},\ldots,X_{\tau_N})$. 
\end{itemize}
When the priority is changing at each step of the game we have the dynamic random priority. The question is if the moment of the assignments, before the arrival of the observation and its presentation to the players of after, has a role. The effective stopping time for player $i$ proposed here assume that the priority is determined before arrival of the observation, and the observation is presented according this order.
\begin{itemize}
\item If the priority is dynamic and random it is defined by the sequence $\Pi=(\Pi_k)_{k=1}^T$. The effective stopping strategy of the player $i$ is following in this case:
\begin{equation}
\tau_i(\vec{p},\Pi)=\inf\{k\geq 1:p_k^i\prod_{j=1}^N(1-p_k^j)\one_{\{j:\Pi_k(j)<\Pi_k(i)\}}=1\},
\end{equation}
It is the stopping time w.r.t. $\tilde{\cF_n^i}=\sigma\{\cF_n^i,\Pi_k,\{(p_k^j)_{k=1,\{j: \Pi_k(j)<\Pi_k(i)\}}^n\}\}$. 
\item Each player stops at different moment. It translates the problem of fixed priority optimization problem to the ordinary stopping game with payoffs 
$$G_i(\tau_1(\vec{p},\Pi),\tau_2(\vec{p},\Pi),\ldots,\tau_N(\vec{p},\Pi))=g_i(X_{\tau_1(\vec{p},\Pi)},X_{\tau_2(\vec{p},\Pi)},\ldots,X_{\tau_N(\vec{p},\Pi)}).
$$ 
\end{itemize}
\vspace{-4ex}
\subsection{Restricted observations of lower priority players.}
\subsubsection{Who has accepted the observation?}
In a sequential decision process taken by the players the consecutive acceptance decision are effectively done by some players. For every stopping time $\tau^i(\vec{p},\Pi)$ the representation by the adapted random sequence $(\delta_k^i)_{k=1}^T$, $i=1,2,\ldots,N$, is given. Let us denote $\gamma_k=\inf\{1\leq i\leq N: \delta_k^i=1\}$ the player who accepted the observation at moment $k$ if any. Similar index  can be defined for the fix deterministic and random priority as for dynamic, deterministic priority as well. 
\subsubsection{Restricted knowledge.}
In the class of such games the natural question which appears is the accessibility of the information. It could be that the accepted observation by the high rang players are hidden for the lower rang players when  has been accepted. However, some information are acquired taking into account the players' behavior.
\begin{itemize}
\item As the result of the decision process players collect information about the states of the process and some of them accept some states. In considered models it was assumed that players are equally informed about the process. Further it will be admitted that the player has access to information according the priority assigned to him. States accepted by others are  not fully accessible to the players which have not seen it before. However, some conjectures are still available assuming rational behavior of the players and it is given by the function $\varphi^i_{\gamma_k}(X_k)$ for the player $i$ when its priority is lower than player's $\gamma_k$th.
\item Effective information available for the player $i$ at moment $k$ can be presented as follows. 
\begin{equation*}
\tilde{X}_k^i=X_k\one_{\{i:i\leq \gamma_k\}}+\varphi_{\gamma_k}^i(X_k)\one_{\{i:i>\gamma_k\}}.
\end{equation*}
\item The player investigation and interaction with other players gives him filtration $\bar{\cF}_n^i=\sigma\{\tilde{\cF}_n^i,\varphi_{\gamma_k}^i(X_k)\}$.
\item Each player stops at different moment. It translates the problem of random priority optimization problem, with restricted access to observation, to the ordinary stopping game with payoffs 
$$G_i(\tau_1(\vec{p},\Pi),\tau_2(\vec{p},\Pi),\ldots,\tau_N(\vec{p},\Pi))=g_i(\tX^i_{\tau_1(\vec{p},\Pi)},\tX^i_{\tau_2(\vec{p},\Pi)},\ldots,\tX^i_{\tau_N(\vec{p},\Pi)}).
$$ 
\end{itemize}

Let us analyze who has accepted the observation? In a sequential decision process taken by the players the consecutive acceptance decision are effectively done by some players. For every stopping time $\tau^i(\vec{p},\Pi)$ the representation by the adapted random sequence $(\delta_k^i)_{k=1}^T$, $i=1,2,\ldots,N$, is given. Let us denote $\gamma_k=\inf\{1\leq i\leq N: \delta_k^i=1\}$ the player who accepted the observation at moment $k$ if any. Similar index  can be defined for the fix deterministic and random priority as for dynamic, deterministic priority as well. 
\section{\label{sec:4}Monotone stopping games with priority}
\subsection{General assumption.} Boundedness assumptions--maximal payoff.
\begin{eqnarray}\vspace{-2.5ex}
\label{MaxPayoff1}
&&\bE(\sup_{\myatop{1\leq j_i\leq T}{i=1,\ldots,N}}G_k(j_1,j_2,\ldots,j_N))< \infty\\
\label{MinPayoff1}\forall_{\myatop{1\leq j_i\leq T}{j\neq i}}&&
\bE(\inf_{1\leq n\leq T}G_k(j_1,\ldots,j_{i-1},n,j_{i+1},\ldots,j_N))>-\infty
\end{eqnarray}
where $k =1,2,\ldots,N$.

In order to assure that each player has a best response to any strategy chosen by other players it is required:
\begin{eqnarray}\vspace{-2.5ex}
\label{assupt3}\forall_{\tau_j\in \cT^j\atop j\neq i} &&\bE[\sup_nG_i(\tau_1,\ldots,\tau_{i-1},n, \tau_{i+1},\ldots,N)|\cF_n^i]^{+}\leq\infty\\
\label{assumpt4}\forall_{\tau_j\in \cT^j\atop j\neq i}&&\limsup_{n\rightarrow T}\bE\left[G_i(\tau_1,\ldots,\tau_{i-1},n,\tau_{i+1},\ldots,N)|\cF_n^i\right]\\
\nonumber&&\leq \bE(G_i(\tau_1,\ldots,\tau_{i-1},T,\tau_{i+1},\ldots,N)|\cF_T^i)\text{ a.e.}.
\end{eqnarray}
For further analyses the conditional expectation of the payoffs for player $i$ should be determine. The sequence $\eta_n=\bE(G_i(\tau_1,\ldots,\tau_{i-1},n,\tau_{i+1},\ldots,\tau_N)|\cF_n)$ is the conditional expected return to player~$i$ if he decide to stop after his $n$th observation and other players uses the stopping rules of their choice. The sequence $\eta_n$ is $\cF_n^i$ adapted and under the boundedness assumption presented above the exists an optimal stopping rule for this sequences for $i=1,2,\ldots,N$. 

\begin{definition}[Regular stopping time] The stopping time $\tau_i\in\cT^i$ is regular with respect to $\tau_1\in\cT^1,\ldots,\tau_{i-1}\in\cT^{i-1},\tau_{i+1}\in\cT^{i+1},\ldots,\tau_N\in\cT^N$ if 
\begin{equation}\label{regstop}
\bE(\eta_{\tau_i}|\cF_n^i)\geq\bE(\eta_n|\cF_n^i)\text{ on $\{\omega:\tau_i>n\}$ for all $n$.}
\end{equation}
Let $\overrightarrow{\tau_{-i}}=(\tau_1,\ldots,\tau_{i-1},\tau_{i+1},\ldots,\tau_n)$.
\end{definition}

Maximal regular best response will be considered. By the results of \cite{mam87:monotone} and \cite{kla73:prop} it can be established:
\begin{lemma}
Under (\ref{MaxPayoff1})-(\ref{assumpt4}) each player has a unique, maximal regular best response $\hat{\tau}_i(\overrightarrow{\tau_{-i}})$ to any vector of stopping times $\overrightarrow{\tau_{-i}}$ chosen by his opponents.
\end{lemma}
This does not immediately imply the Nash equilibrium existence.
\subsection{Monotone structure of best responses.}
The incremental benefit to player should be analyzed. It is assumed that increments of payoffs have the following properties. Let us consider the following increments of payoffs.
\begin{eqnarray*}
\forall_{m<T\atop 1\leq k\leq T}\Delta_m^i(k,\overrightarrow{j_{-i}})&=&G_i(j_1,\ldots,j_{i-1},m+k,j_{i+1},\ldots,j_N)\\
&&-G_i(j_1,\ldots,j_{i-1},m,j_{i+1},\ldots,j_N) 
\end{eqnarray*}
\begin{description}
\item[ND] Let us assume that  $\Delta_m^i(k,\vec{j}_{-i})$ is nondecreasing in $\overrightarrow{j_{-i}}$;
\item[NI] Let us assume that  $\Delta_m^i(k,\vec{j}_{-i})$ is nonincreasing in $\overrightarrow{j_{-i}}$;
\end{description}

\begin{lemma}
If (\ref{MaxPayoff1})-(\ref{assumpt4}) and condition ND are fulfilled and $\sigma\in\cT^i$ is regular with respect of $\overrightarrow{\tau_{-i}}_1\in \cT^{-i}$ then it is also regular with respect to any $\overrightarrow{\tau_{-i}}_2\in \cT^{-i}$ such that $\overrightarrow{\tau_{-i}}_1\preceq\overrightarrow{\tau_{-i}}_2$ a.e.. (Under NI $\overrightarrow{\tau_{-i}}_2\preceq\overrightarrow{\tau_{-i}}_1$ a.e.)    
\end{lemma}
\begin{lemma}
Let $\overrightarrow{\tau_{-i}}_k\in\cT^{-i}$, $k=1,2$, such that $\overrightarrow{\tau_{-i}}_1\preceq\overrightarrow{\tau_{-i}}_2$ a.e. and ND is fulfilled then the best response $\hat{\sigma}(\overrightarrow{\tau_{-i}}_1)\preceq\hat(\sigma)(\overrightarrow{\tau_{-i}}_2)$  a.e.. (Under NI $\hat{\sigma}(\overrightarrow{\tau_{-i}}_1)\succeq\hat{\sigma}(\overrightarrow{\tau_{-i}}_2)$ a.e..)
\end{lemma}
\subsubsection{Tarski's fixed point theorem.}
The fixed point theorem which will be helpful for proving the existence of the equilibrium is obtained for the complete lattices and an isotone functions. We consider the partial order of random variables: $\tau\preceq\sigma$ iff $\tau\leq\sigma$ \emph{a.e.}. The operations of supremum of random variables and infimum of random variables are inner operations in $\cT^i$. If the essential supremum is considered we have also for every subset $\cA\subset \cS$ that $\vee\cA\in\cS$ and $\wedge\cA\in\cS$.
\begin{lemma}[Stopping set is a complete lattice]
The partially ordered sets $\cT^i$ with order $\preceq$ and operations essential supremum $\vee$ and essential infimum  $\wedge$ defined in it are complete lattices.
\end{lemma}
\begin{definition}[Isotone function]
Let $\cS$ be lattice. $f$ is isotone function from $\cS$ into $\cS$ if for $\tau,\sigma\in\cS$ such that $\tau\preceq\sigma$ implies $f(\tau)\leq f(\sigma)$.
\end{definition}
\begin{theorem}[\cite{tar55:fixed}]
If $\cS$ is a complete lattice and if $f$ is an isotone function from $\cS$ into $\cS$, then $f$ has a fixed point.
\end{theorem}

\subsection{Main result.}
Monotonicity of increments with integrability of payoff functions guarantee existence of Nash equilibrium in stopping game with various models of priority (rule of assignments) based on the theorem.  
\begin{theorem}[\cite{mam87:monotone}]
Suppose that assumptions (\ref{MaxPayoff1})-(\ref{assumpt4}) with ND or NI holds. Then there is a Nash equilibrium pair of stopping times. There is an vector of stopping times $\vec{\tau^\star}$ which forms an equilibrium such that $\tau_i^\star=\hat{\sigma}(\overrightarrow{\tau_{-i}^\star})$, $i=1,2,\ldots,N$.
\end{theorem}

\section{\label{sec:5}Conclusion.} Based on the consideration of the paper we know that the various priority approach model in the multiple choice problem can be transformed to the multiperson antagonistic game with the equilibrium point as the rational treatment. The equilibrium point in all these problems exist. The construction of them need individual tretment and it is not solved in general yet.  

The presented decision model can be found with a slightly different interpretation, namely games with the arbitration procedure. Details can be found e.g. in the works of \citeauthor{Sak1984:Arbitration}(\citeyear{Sak1984:Arbitration}) and \citeauthor{MazSakZab2002:Arbitration}(\citeyear{MazSakZab2002:Arbitration}).\footnote{See also \citeauthor{BraMer1992:Arbitration}(\citeyear{BraMer1992:Arbitration}) and \citeauthor{Cha1981:Arbitration}(\citeyear{Cha1981:Arbitration}) for details concerning arbitration procedure.}

The close to the models are some multivariate stopping problem with cooperation. In cooperative stopping games the players have to use the decision suggested by coordinator of the decision process (cf. \citeauthor{asssam98:multi}(\citeyear{asssam98:multi}), \citeauthor{gli04:coop}(\citeyear{gli04:coop}). 

In \citeauthor{kuryasnak80:multi}(\citeyear{kuryasnak80:multi}) the idea of voting stopping rules has been proposed. The game defined on the sequence of iid random vectors has been defined with the concept of the Nash equilibrium as the solution. There are generalization of the results obtained by \citeauthor{szayas95:voting}(\citeyear{szayas95:voting}).  Conditions for a unique equilibrium among stationary threshold strategies in such games are given by  \citeauthor{fer04:voting}(\citeyear{fer04:voting}).     

\section*{Acknowledgments} The authors' thanks go to many colleagues taking part in discussion of the topics presented in the paper.


\begin{thebibliography}{44}
\providecommand{\natexlab}[1]{#1}
\providecommand{\url}[1]{\texttt{#1}}
\expandafter\ifx\csname urlstyle\endcsname\relax
  \providecommand{\doi}[1]{doi: #1}\else
  \providecommand{\doi}{doi: \begingroup \urlstyle{rm}\Url}\fi

\bibitem[Assaf and Samuel-Cahn(1998)]{asssam98:coop}
D.~Assaf and E.~Samuel-Cahn.
\newblock {Optimal cooperative stopping rules for maximization of the product
  of the expected stopped values.}
\newblock \emph{Stat. Probab. Lett.}, 38\penalty0 (1):\penalty0 89--99, 1998.
\newblock \doi{10.1016/S0167-7152(97)00158-2}.
\newblock \ZBL{0912.60061}.

\bibitem[{Assaf} and {Samuel-Cahn}(1998)]{asssam98:multi}
D.~{Assaf} and E.~{Samuel-Cahn}.
\newblock {Optimal multivariate stopping rules.}
\newblock \emph{{J. Appl. Probab.}}, 35\penalty0 (3):\penalty0 693--706, 1998.
\newblock ISSN 0021-9002; 1475-6072/e.
\newblock \doi{10.1239/jap/1032265217}.
\newblock \ZBL{0937.60040}.

\bibitem[{Brams} and {Merrill, III}(1992)]{BraMer1992:Arbitration}
S.~J. {Brams} and S.~{Merrill, III}.
\newblock
  \href{http://control.ibspan.waw.pl:3000/contents/export?filename=1992-1-09_brams_merrill.pdf}{Arbitration
  procedures with the possibility of compromise}.
\newblock In J.~Stefański, editor, \emph{Bargaining and arbitration in
  conflicts}, volume 21(1) of \emph{Control and Cybernetics}, pages 131--149.
  Systems Research Institute of the Polish Academy of Sciences, 1992.
\newblock ISSN {0324-8569}; \MR{1218889}.

\bibitem[Chatterjee(1981)]{Cha1981:Arbitration}
K.~Chatterjee.
\newblock Comparison of arbitration procedures: models with complete and
  incomplete information.
\newblock \emph{IEEE Trans. Systems Man Cybernet.}, 11\penalty0 (2):\penalty0
  101--109, 1981.
\newblock ISSN 0018-9472.
\newblock \doi{10.1109/TSMC.1981.4308635}.
\newblock \MR{611435}.

\bibitem[{Diecidue} and {van de Ven}(2008)]{DieVen2008:Aspiration}
E.~{Diecidue} and J.~{van de Ven}.
\newblock Aspiration level, probability of success and failure, and expected
  utility.
\newblock \emph{Internat. Econom. Rev.}, 49\penalty0 (2):\penalty0 683--700,
  2008.
\newblock ISSN 0020-6598.
\newblock \doi{10.1111/j.1468-2354.2008.00494.x}.
\newblock \MR{2404450}.

\bibitem[{Dorobantu} et~al.(2009){Dorobantu}, {Mancino}, and
  {Pontier}]{DorManPon2009:OS}
D.~{Dorobantu}, M.~E. {Mancino}, and M.~{Pontier}.
\newblock Optimal strategies in a risky debt context.
\newblock \emph{Stochastics}, 81\penalty0 (3-4):\penalty0 269--277, 2009.
\newblock ISSN 1744-2508.
\newblock \doi{10.1080/17442500902917433}.
\newblock \MR{2549487}.

\bibitem[{Dynkin}(1969)]{dyn69:game}
E.~{Dynkin}.
\newblock {Game variant of a problem on optimal stopping.}
\newblock \emph{{Sov. Math., Dokl.}}, 10:\penalty0 270--274, 1969.
\newblock ISSN 0197-6788.
\newblock Translation from Dokl. Akad. Nauk SSSR 185, 16-19 (1969).
  \ZBL{0186.25304}.

\bibitem[Enns and Ferenstein(1987)]{EnnFer1987:Priorities}
E.~G. Enns and E.~Z. Ferenstein.
\newblock On a multiperson time-sequential game with priorities.
\newblock \emph{Sequential Anal.}, 6\penalty0 (3):\penalty0 239--256, 1987.
\newblock ISSN 0747-4946.
\newblock \doi{10.1080/07474948708836129}.
\newblock \MR{918908}.

\bibitem[{Etro}(2013)]{Etro2013:ReviewStackelberg}
F.~{Etro}.
\newblock {Book review of: Heinrich von Stackelberg, Market structure and
  equilibrium.}
\newblock \emph{{J. Econ.}}, 109\penalty0 (1):\penalty0 89--92, 2013.
\newblock ISSN 0931-8658; 1617-7134/e.
\newblock \doi{10.1007/s00712-013-0341-9}.

\bibitem[Feng and Xiao(2000)]{FenXia2000:Revenue}
Y.~Feng and B.~Xiao.
\newblock Revenue management with two market segments and reserved capacity for
  priority customers.
\newblock \emph{Adv. in Appl. Probab.}, 32\penalty0 (3):\penalty0 800--823,
  2000.
\newblock ISSN 0001-8678.
\newblock \doi{10.1239/aap/1013540245}.
\newblock \MR{1788096}.

\bibitem[Ferenstein(1992)]{Fer1992:2person}
E.~Z. Ferenstein.
\newblock Two-person non-zero-sum sequential games with priorities.
\newblock In \emph{Strategies for sequential search and selection in real time
  ({A}mherst, {MA}, 1990)}, volume 125 of \emph{Contemp. Math.}, pages
  119--133. Amer. Math. Soc., Providence, RI, 1992.
\newblock \doi{10.1090/conm/125/1160615}.
\newblock \MR{1160615}.

\bibitem[Ferguson(2005)]{fer04:voting}
T.~Ferguson.
\newblock Selection by committee.
\newblock In K.~S. A.S.~Nowak, editor, \emph{Advances in dynamic games:
  applications to economics, finance, optimization, and stochastic control},
  volume~7 of \emph{Annals of the International Society of Dynamic Games},
  pages 203--209, Boston, 2005. Birkh\"aser.
\newblock \ZBL{1123.91003}.

\bibitem[Ferguson(2016)]{Fer2016:SGSakaguchi}
T.~S. Ferguson.
\newblock The sum-the-odds theorem with application to a stopping game of
  {S}akaguchi.
\newblock \emph{Math. Appl. (Warsaw)}, 44\penalty0 (1):\penalty0 45--61, 2016.
\newblock ISSN 1730-2668.
\newblock \doi{10.14708/ma.v44i1.1192}.
\newblock \MR{3557090}.

\bibitem[Fushimi(1981)]{fus1981:competitive}
M.~Fushimi.
\newblock The secretary problem in a competitive situation.
\newblock \emph{J. Oper. Res. Soc. Jap.}, 24:\penalty0 350--358, 1981.
\newblock \ZBL{0482.90090}.

\bibitem[Glickman(2004)]{gli04:coop}
H.~Glickman.
\newblock {Cooperative stopping rules in multivariate problems.}
\newblock \emph{Sequential Anal.}, 23\penalty0 (3):\penalty0 427--449, 2004.

\bibitem[Haggstrom(1967)]{hagg1967:2stop}
G.~Haggstrom.
\newblock Optimal sequential procedures when more then one stop is required.
\newblock \emph{Ann. Math. Statist.}, 38:\penalty0 1618--1626, 1967.

\bibitem[{Jaśkiewicz} and {Nowak}(2016)]{JasNow2016:Non-zero}
A.~{Jaśkiewicz} and A.~{Nowak}.
\newblock Non-zero-sum stochastic games.
\newblock In T.~Ba\c{s}sar and G.~Zaccour, editors, \emph{Handbook of Dynamic
  Game Theory}, page 64pp. Birkhäuser and Springer International Publishing AG
  of Springer Nature, Bassel, 2016.
\newblock \doi{10.1007/978-3-319-27335-8\_33-1}.

\bibitem[Kifer(1969)]{Kif1969:Games}
J.~I. Kifer.
\newblock Optimal behavior in games with an infinite sequence of moves.
\newblock \emph{Teor. Verojatnost. i Primenen.}, 14:\penalty0 284--291, 1969.
\newblock ISSN 0040-361x.

\bibitem[{Klass}(1973)]{kla73:prop}
M.~J. {Klass}.
\newblock {Properties of optimal extended-valued stopping rules for S$_n$/n.}
\newblock \emph{{Ann. Probab.}}, 1:\penalty0 719--757, 1973.
\newblock ISSN 0091-1798; 2168-894X/e.
\newblock \doi{10.1214/aop/1176996843}.
\newblock \ZBL{0281.62085}.

\bibitem[Krasnosielska-Kobos and Ferenstein(2013)]{KraFer2013:Elfving}
A.~Krasnosielska-Kobos and E.~Ferenstein.
\newblock Construction of {N}ash equilibrium in a game version of {E}lfving's
  multiple stopping problem.
\newblock \emph{Dyn. Games Appl.}, 3\penalty0 (2):\penalty0 220--235, 2013.
\newblock ISSN 2153-0785.
\newblock \doi{10.1007/s13235-012-0070-7}.

\bibitem[Kurano et~al.(1980)Kurano, Yasuda, and Nakagami]{kuryasnak80:multi}
M.~Kurano, M.~Yasuda, and J.~Nakagami.
\newblock Multi-variate stopping problem with a majority rule.
\newblock \emph{J. Oper. Res. Soc. Jap.}, 23:\penalty0 205--223, 1980.
\newblock \ZBL{0439.90105}.

\bibitem[Mamer(1987)]{mam87:monotone}
J.~W. Mamer.
\newblock {Monotone stopping games.}
\newblock \emph{{J. Appl. Probab.}}, 24:\penalty0 386--401, 1987.
\newblock ISSN 0021-9002; 1475-6072/e.
\newblock \doi{10.2307/3214263}.
\newblock \ZBL{0617.90092}.

\bibitem[Mazalov et~al.(2002)Mazalov, Sakaguchi, and
  Zabelin]{MazSakZab2002:Arbitration}
V.~V. Mazalov, M.~Sakaguchi, and A.~A. Zabelin.
\newblock Multistage arbitration game with random offers.
\newblock \emph{Int. J. Math. Game Theory Algebra}, 12\penalty0 (5):\penalty0
  409--417, 2002.
\newblock ISSN 1060-9881.
\newblock \MR{1951832}; \ZBL{1076.91502}.

\bibitem[{McKean, Jr., H. P.}(1965)]{McK1965:Free}
{McKean, Jr., H. P.}
\newblock {Appendix: A free boundary problem for the heat equation arising from
  a problem of mathematical economics.}
\newblock \emph{{Ind. Manag. Rev.}}, 6:\penalty0 32--39, 1965.

\bibitem[Neumann et~al.(1994)Neumann, Porosiński, and
  Szajowski]{neuporsza94:note}
P.~Neumann, Z.~Porosiński, and K.~Szajowski.
\newblock {A} note on two person full-information best chice problems with
  imperfect observation.
\newblock In \emph{Operations Research. Extended Abstracts of the 18th
  Symposium on Operations Research (GMOOR), Cologne, sept. 1-3, 1993.}, pages
  355--358, Heidelberg, 1994. GMOOR, Phisica-Verlag.

\bibitem[Neumann et~al.(2002)Neumann, Ramsey, and
  Szajowski]{neuramsza02:Randomized}
P.~Neumann, D.~Ramsey, and K.~Szajowski.
\newblock Randomized stopping times in {D}ynkin games.
\newblock \emph{ZAMM Z. Angew. Math. Mech.}, 82\penalty0 (11-12):\penalty0
  811--819, 2002.
\newblock ISSN 0044-2267.
\newblock \doi{10.1002/1521-4001(200211)82:11/12<811::AID-ZAMM811>3.0.CO;2-P}.
\newblock 4th GAMM-Workshop ``Stochastic Models and Control Theory''
  (Lutherstadt Wittenberg, 2001).

\bibitem[Nowak and Szajowski(1998)]{nowsza94:stochastic}
A.~Nowak and K.~Szajowski.
\newblock {N}onzero-sum stochastic games.
\newblock In T.~P. M.~Bardi, T.E.S.~Raghavan, editor, \emph{Stochastic and
  Differential Games. Theory and Numerical Methods}, Annals of the
  International Society of Dynamic Games, pages 297--342, Boston, 1998.
  Birkh\"{a}ser.
\newblock \MR{200d:91021}; \ZBL{0940.91014}.

\bibitem[Radzik and Szajowski(1988)]{RadSza1988:sequential}
T.~Radzik and K.~Szajowski.
\newblock On some sequential game.
\newblock \emph{Pure Appl. Math. Sci.}, 28\penalty0 (1-2):\penalty0 51--63,
  1988.
\newblock ISSN 0379-3168.

\bibitem[Ramsey and Cierpia\l(2009)]{RamCie2009:Cooperative}
D.~Ramsey and D.~Cierpia\l.
\newblock Cooperative strategies in stopping games.
\newblock In \emph{Advances in dynamic games and their applications}, volume~10
  of \emph{Ann. Internat. Soc. Dynam. Games}, pages 415--430. Birkh\"auser
  Boston, Inc., Boston, MA, 2009.

\bibitem[Ramsey and Szajowski(2008)]{ramsza08:correlated}
D.~M. Ramsey and K.~Szajowski.
\newblock {Selection of a correlated equilibrium in Markov stopping games.}
\newblock \emph{Eur. J. Oper. Res.}, 184\penalty0 (1):\penalty0 185--206, 2008.
\newblock \doi{10.1016/j.ejor.2006.10.050}.

\bibitem[Ravindran and Szajowski(1992)]{RavSza1992:Dynkin}
G.~Ravindran and K.~Szajowski.
\newblock Nonzero sum game with priority as {D}ynkin's game.
\newblock \emph{Math. Japon.}, 37\penalty0 (3):\penalty0 401--413, 1992.
\newblock ISSN 0025-5513.

\bibitem[Sakaguchi(1984)]{Sak1984:Arbitration}
M.~Sakaguchi.
\newblock A time-sequential game related to an arbitration procedure.
\newblock \emph{Math. Japon.}, 29\penalty0 (3):\penalty0 491--502, 1984.
\newblock ISSN 0025-5513.

\bibitem[{Sakaguchi}(1995)]{sak95:review}
M.~{Sakaguchi}.
\newblock {Optimal stopping games -- a review.}
\newblock \emph{{Math. Japon.}}, 42\penalty0 (2):\penalty0 343--351, 1995.
\newblock ISSN 0025-5513.
\newblock {Correction to ``Optimal stopping games -- a review''.
  \ZBL{0879.60044}} \ZBL{0865.60035}.

\bibitem[{Shmaya} and {Solan}(2004)]{shmsol04:nonzero}
E.~{Shmaya} and E.~{Solan}.
\newblock {Two-player nonzero-sum stopping games in discrete time.}
\newblock \emph{{Ann. Probab.}}, 32\penalty0 (3B):\penalty0 2733--2764, 2004.
\newblock ISSN 0091-1798; 2168-894X/e.
\newblock \doi{10.1214/009117904000000162}.
\newblock \ZBL{1079.60045}.

\bibitem[Szajowski(1992)]{Sza1992:Priority}
K.~Szajowski.
\newblock On non-zero sum game with priority in the secretary problem.
\newblock \emph{Math. Japonica}, 37\penalty0 (3):\penalty0 415--426, 1992.

\bibitem[Szajowski(1993{\natexlab{a}})]{Sza1993:2S2DM}
K.~Szajowski.
\newblock Double stopping by two decision-makers.
\newblock \emph{Adv. in Appl. Probab.}, 25\penalty0 (2):\penalty0 438--452,
  1993{\natexlab{a}}.
\newblock ISSN 0001-8678.
\newblock \doi{10.2307/1427661}.

\bibitem[Szajowski(1993{\natexlab{b}})]{Sza1993:ZOR}
K.~Szajowski.
\newblock {Markov} stopping games with random priority.
\newblock \emph{Zeitschrift f\"ur Operations Research}, 37\penalty0
  (3):\penalty0 69--84, 1993{\natexlab{b}}.

\bibitem[Szajowski(1995)]{Sza1995:SIAM}
K.~Szajowski.
\newblock Optimal stopping of a discrete {M}arkov process by two decision
  makers.
\newblock \emph{SIAM J. Control Optim.}, 33\penalty0 (5):\penalty0 1392--1410,
  1995.
\newblock ISSN 0363-0129.
\newblock \doi{10.1137/S0363012993246877}.

\bibitem[{Szajowski} and {Yasuda}(1997)]{szayas95:voting}
K.~{Szajowski} and M.~{Yasuda}.
\newblock {Voting procedure on stopping games of Markov chain.}
\newblock In S.~O. Anthony H.~Christer and L.~C. Thomas, editors,
  \emph{UK-Japanese Research Workshop on Stochastic Modelling in Innovative
  Manufacturing, July 21-22, 1995}, volume 445 of \emph{Lecture Notes in
  Economics and Mathematical Systems}, pages 68--80. Moller Centre, Churchill
  College, Univ. Cambridge, UK, Springer, 1997.
\newblock ISBN 3-540-61768-X/pbk.
\newblock \MR{98a:90159}; \ZBL{0878.90112}.

\bibitem[{Tarski}(1955)]{tar55:fixed}
A.~{Tarski}.
\newblock {A lattice-theoretical fixpoint theorem and its applications.}
\newblock \emph{{Pac. J. Math.}}, 5:\penalty0 285--309, 1955.
\newblock ISSN 0030-8730.
\newblock \doi{10.2140/pjm.1955.5.285}.
\newblock \ZBL{0064.26004}.

\bibitem[{von Stackelberg}(2011)]{Sta2011:enMarket}
H.~{von Stackelberg}.
\newblock \emph{{Market structure and equilibrium. Translated from the German
  by Damian Bazin, Lynn Urch and Rowland Hill.}}
\newblock Berlin: Springer, 2011.
\newblock ISBN 978-3-642-12585-0/hbk; 978-3-642-12586-7/ebook.
\newblock \doi{10.1007/978-3-642-12586-7}.

\bibitem[von Stackelberg(1934)]{Sta1934:Marktform}
H.~F. von Stackelberg.
\newblock \emph{{Marktform und Gleichgewicht}}.
\newblock Springer, Wien and Berlin, 1934.

\bibitem[von Stackelberg(2009)]{Sta1932:Kostentheorie}
H.~F. von Stackelberg.
\newblock \emph{{Grundlagen einer reinen Kostentheorie}}.
\newblock {Meilensteine Nationalokonomie}. Springer-Verlag Gmbh, Berlin, 2009.
\newblock Originally published monograph. Reprint of the 1st Ed. Wien, Verlag
  von Julius Springer, 1932.
  \href{http://www.springerlink.com/content/978-3-540-85271-1}{Read on line}.

\bibitem[{Yasuda}(1985)]{yas85:randomize}
M.~{Yasuda}.
\newblock {On a randomized strategy in Neveu's stopping problem.}
\newblock \emph{{Stochastic Processes Appl.}}, 21:\penalty0 159--166, 1985.
\newblock ISSN 0304-4149.
\newblock \doi{10.1016/0304-4149(85)90384-9}.
\newblock \ZBL{0601.60039}.

\end{thebibliography}
\end{document}